\documentclass{bull-l}

\usepackage {hanging}
\usepackage {url}

\begin{document}

\title{}

\author{Joseph F. Grcar}
\address{6059 Castlebrook Drive; Castro Valley, CA 94552-1645 USA}
\email{jfgrcar@gmail.com}


\maketitle

\noindent \hangindent=\parindent \hangafter=1 
\textit {The Chinese Roots of Linear Algebra\/}, by Roger Hart, The Johns Hopkins University Press, Baltimore, 2011, xiii + 286 pp., \$65.00, ISBN 8-8018-9755-6

\bigskip

A recurring theme in the history of mathematics is the cultural setting of the work. Believe it or not, social context is controversial \cite {Rowe1996}. On one side, context must be important because mathematics rests on human experience \cite {vonNeumann1947b-short}; on the other hand, society cannot matter because mathematics is pure thought. The latter view is reinforced by what specialists learn of the past through their work: a ``heritage'' \cite {GrattanGuinness2004} of famous ideas  attributed to the past but couched in the latest terminology, which makes the ideas appear to be immutable. \textit {The Chinese Roots of Linear Algebra\/} by Roger Hart will stimulate some lively discussion because many books about the heritage of mathematics barely mention China.

Several types of calculations were performed in ancient China that continued to be practiced in eastern Asia through the 19th century. The calculations are described in historical studies, e.g.\ \cite {Dauben2007, Li1987, Martzloff1997-BAMS}, and they are occasionally mentioned in mathematics texts, e.g.\ \cite [p.\ 1] {Meyer2000}. Some ancient calculations resemble those of the present day, including the ``Gaussian elimination'' of American precalculus algebra textbooks. This similarity in different cultures raises the question whether they share an immutable mathematical concept and what it might be (Hart, p.\ 191), in that one development involves symbolic algebra, and the other did not.

The ancient calculations were made with implements called counting rods (the abacus being a comparatively recent invention \cite [chap.\ 6] {Li1987}). The rods were short sticks that could be arranged to represent digits in a decimal notation for positive and negative integers. Elaborate calculations were made by placing the rods inside squares on a counting ``board'' or ``table.'' No ancient tables survive \cite [chap.\ 13] {Martzloff1997-BAMS}, but by one surmise, they were any flat surface, perhaps covered by a sheet of cloth ruled into squares. The counting table was literally an ancient ``spreadsheet'' for manual computing in which numbers could be entered and changed as a calculation progressed. 

Of sources for mathematics that survive from ancient China, the most comprehensive is the \textit {Nine Chapters of the Mathematical Art\/}, which was compiled anonymously about two millennia ago. Like many ancient texts, the \textit {Nine Chapters\/} has been reassembled from fragmentary copies that were made long after the reputed origin. Some reconstructions can be read in translation, e.g.\ \cite {Chemla2005, Shen1999}. Each of the nine chapters treats a different type of calculation through didactic examples. Most of the lessons are \textit {l'art pour l'art\/}, in that uses for the knowledge are not identified. (Other texts describe administrative and astronomical mathematics \cite [p.\ 187] {Dauben2007}.) The problems were supplemented by explanations in derivative works. The earliest of these ``commentaries'' was written by Liu Hui when the \textit {Nine Chapters\/} were about two hundred years old. 

Roger Hart's focus is Chapter Eight of the \textit {Nine\/}. These are the problems interpretable in modern algebra as simultaneous linear equations. Among Hart's theses is that Chapter Eight represents the practices of ``illiterate adepts,'' whose work was recorded by ``aspiring literati,'' who then presented their writings to the royal court (p.\ 1). This view has been stated elsewhere for more recent periods \cite [chap.\ 9] {Martzloff1997-BAMS}. It would be fascinating to know more about academic career paths in historic China. The thesis does have parallels in our society. Today, ``a critical point in the career of the young mathematician'' is being recognized by the central government \cite {Brennan2007}. Moreover, before electronic computers, calculations were made by anonymous hand computers \cite {Grier2005}, and now, computer programmers work in similar obscurity.

\renewcommand {\a} [3] {\alpha^{(#1)}_{#2,#3}}
\renewcommand {\t} [1] {A^{(#1)}}

The linear calculations described in Chapter Eight were as follows. So as to avoid circumlocutions, note the blanket caveat that this explanation interprets ancient mathematics in modern terms. A system of $n$ linear equations with the same quantity of unknowns $x_j$, 
\begin {displaymath}
\sum_{j=1}^n a_{i,j} x_j = b_i \quad \mbox {for $i = 1, \dots, n$},
\end {displaymath}
was represented by placing the numbers for each equation in a vertical column of the counting table. Rotating columns by $\pi / 2$ to give rows makes the modern tableau, for $n = 3$,
\begin {displaymath}
\t 1 = 
\left[ \begin {array} {cccc} 
a_{1,1}& a_{1,2}& a_{1,3}& b_1\\
a_{2,1}& a_{2,2}& a_{2,3}& b_2\\
a_{3,1}& a_{3,2}& a_{3,3}& b_3
\end {array} \right] .
\end {displaymath}
The ancient calculation amounts to multiplying the tableau on the left by the matrix 
\begin {displaymath}
\left[ \begin {array} {cccc} 
1& 0& 0\\
- a_{2,1}& a_{1,1}& 0\\
- a_{3,1}& 0& a_{1,1}
\end {array} \right] 
\end {displaymath}
which creates a second tableau $\t 2$ of integer entries with $0$ in the first column below the diagonal. I will use the uniform notation $\smash {\a k i j}$ for the entries of the $k$-th tableau. Upon $n-1$ such steps the counting-board cum tableau is in the form algebraists call ``row echelon,'' and numerical analysts call ``upper triangular,'' 
\begin {displaymath}
\renewcommand {\arraystretch} {1.5}
\t 3 =
\left[ \begin {array} {ccccc} 
\a 3 1 1& \a 3 1 2& \a 3 1 3& \a 3 1 4\\
0& \a 3 2 2& \a 3 2 3& \a 3 2 4\\
0& 0& \a 3 3 3& \a 3 3 4\\
\end {array} \right] .
\end {displaymath}
Converting back to symbolic algebra, this form of the tableau corresponds to recurrence formulas from which the unknowns can be evaluated in reverse order,
\begin {displaymath}
\sum_{j=i}^n \a n i j x_j = \a n i {n+1} \quad \mbox {for $i = n, \dots, 1$}.
\end {displaymath}

It is remarkable that Gaussian elimination appears in modern precalculus textbooks with the same visual style of ancient China. The much later, earliest example yet found in Europe is by Jean Borrel in 1560 \cite [pp.\ 190+] {Buteo1560-Notices}. He preserved integers by the same reduction process, and then he evaluated the unknowns using the recurrences in a process now called ``back substitution.'' He used these processes to solve three equations and unknowns. Although the calculation may look like what we do, the concept of ``equation'' was not fully formed when Borrel wrote \cite {Heeffer2011}. For four equations and unknowns, Borrel began with symbolic elimination, but then he reverted to verbal argumentation without finishing the reduction process. 

Because back substitution presumes symbolic algebra, Hart suggests the conclusion of the ancient counting-table calculation has been misunderstood (p.\ 93). He believes it continued the same sort of reduction process to bring the tableau into diagonal form, with the benefit of preserving integers even when the solution has fractions. Columns from right to left can be cleared above the diagonal, and the magnitudes of the remaining numbers can be moderated with divisions. In one example (the oft-exhibited Problem 1 of Chapter 8), the first backward step amounted to multiplying $\t n$ on the left by the matrix
\begin {displaymath}
\renewcommand {\arraystretch} {1.5}
\left[ \begin {array} {ccccc} 
\a 3 1 1& 0& \;0\;\\
0& \a 3 2 2& 0\\
0& 0& 1
\end {array} \right]^{-1} 
\left[ \begin {array} {ccccc} 
\a 3 3 3& 0& - \a 3 2 4\\
0& \a 3 3 3& - \a 3 3 4\\
0& 0& 1
\end {array} \right] 
\end {displaymath}
which creates the next tableau, $\t 4$, with $0$ everywhere off the diagonal in column $3$. Hart shows the division preserves integers. After $n-1$ of these backward steps, the equations have been placed in diagonal form with $\smash {\a n n n}$ replicated on the diagonal,
\begin {displaymath}
\renewcommand {\arraystretch} {1.5}
\t 5 =
\left[ \begin {array} {ccccc} 
\a 3 3 3& 0& 0& \a 5 1 4\\
0& \a 3 3 3& 0& \a 4 2 4\\
0& 0& \a 3 3 3& \a 3 3 4\\
\end {array} \right] .
\end {displaymath}
Hart's thesis about deferring the formation of fractions is compelling. Whether the divisions were performed exactly as he hypothesizes seems to rest on the very few complete examples of the \textit {Nine Chapters\/}. 

Among the problems in Chapter Eight are several of a special form whose tableau representation is, for $n = 4$,
\begin {displaymath}
\left[ \begin {array} {cccc|c} 
a_1& 1& 0& 0& b\\
0& a_2& 1& 0& b\\
0& 0& a_3& 1& b\\
1& 0& 0& a_4& b
\end {array} \right] 
=
\left[ \begin {array} {c | c}
A& y
\end {array} \right] ,
\end {displaymath}
where $A$ and $y$ are the indicated matrix and column vector. These problems have from three to five equations in the \textit {Nine Chapters\/}, and more in later examples. If the reduction phase is done in the straightforward way, with integers, then the final diagonal entry of the row echelon form is $\smash {\a n n n} = \det (A)$. One such problem, the ``well problem,'' has $b$ as an unknown, giving more unknowns than equations. This one example has often been cited to mean that indeterminate problems were understood by the author(s) of the \textit {Nine Chapters\/}. Hart argues to the contrary that the solution in Chapter Eight merely posits the value $b = \det (A)$. A survey of extant treatises (which fills the entire appendix B) reveals indeterminacy might not have been explicitly treated until the 17th century (p.\ 122). 

Hart further suggests the ``determinantal-style'' solutions of the special problems may have affected European work on determinants, and more broadly on linear algebra, through Gottfried Wilhelm Leibniz (pp.\ 26, 189). It is left for future work to find Chinese influences in the mathematical writings of Leibniz. Whatever they may be, Leibniz did not contribute to Gaussian elimination \cite {Grcar2011c, Grcar2011e}, which appeared without determinants in the work of several European mathematicians such as Borrel. Isaac Newton \cite [60--62] {Newton1720-Notices} began a tradition of creating elimination lessons for instruction in symbolic algebra, which culminated in a standardized textbook presentation by the end of the 18th century. This lesson is the ``Gaussian elimination'' of precalculus algebra textbooks. 

In summary, \textit {The Chinese Roots of Linear Algebra\/} chronicles the linear problems of ancient China in the \textit {Nine Chapters\/}, and supplies new insights about their solution. What remains to investigate is whether Chapter Eight of the \textit {Nine\/} influenced modern linear algebra. Are the \textit {Nine Chapters\/} a ``root,'' or are they a separate development, and either way, are they not part of our mathematical heritage? Roger Hart's provocative book deserves to be in every college and university collection. The author's own studies, and his assessment of other scholarly work, will be starting points for innumerable term papers. Beyond that, the specific topic of the title makes engaging reading.

\section* {Appendix}

Perhaps it should be mentioned: the best that can be done with integer elimination is Chi\`o's method \cite {Chio1853, Fuller1975}. The reduction to row echelon form of a tableau $\smash {A^{(1)}} = \smash {\big [ \a 1 i j \big ]}$ with $n$ rows (indexed by the first subscript), and as many or more columns (indexed by the second), can be expressed by a formula for the changed entries of the successive tableaux (indexed by the superscript),
\begin {displaymath}
\a {k+1} i j = {1 \over \Delta (k)} \, \det \left[ \begin {array} {cc} \a k k k& \a k k j\\ \noalign {\medskip} \a k i k& \a k i j\end {array} \right] \quad \mbox {for} \left\{ \, \parbox {8em} {\strut $k = 1, \ldots, n - 1$,\\ \mbox {then $j \ge k$,}\\ \strut \mbox {and then $i > k$,}} \right.
\end {displaymath}
where $\Delta (k)$ is to be chosen. (Note the case $j = k$ introduces $0$ into the columns.) The \textit {Nine Chapters\/} has $\Delta (k) = 1$, a conventional form of Gaussian elimination has $\Delta (k) = \smash {\a k k k}$, and Chi\`o has $\Delta (k) = \smash {\a {k-1} {k-1} {k-1}}$ or $1$ for $k = 1$. For Chi\`o's choice, it can be shown: (1) the tableaux entries remain integers, (2) for $i, j \ge k$, the degree of $\smash {\a k i j}$ as a polynomial in the entries of the initial tableau is $k$, which is the smallest that can be expected for any calculation, and (3) the entry $\smash {\a k k k}$ is the leading principle minor of order $k$ for the initial tableau. 

Many authors have reinvented Chi\`o's method for evaluating determinants, including Charles Dodgson (Lewis Carroll). The method is also attributed to Bareiss \cite {Bareiss1968} whose novel contribution was to accelerate the process with higher-order determinants. Chi\`o's method of course applies to elimination in integral domains in which form it may be found in commutative algebra and complexity theory.

If $i$ is permitted to range over $1, \ldots, n$ except $k$ then the $n$-th tableau is diagonal, and the reduction is called Gauss-Jordan elimination \cite {Althoen1987, Katz1988} as opposed to Gaussian elimination. The Jordan version requires more arithmetic operations to obtain a diagonal form than reduction to row echelon form followed by a backward phase.

\section* {Acknowledgments}

I thank William Kahan and Jonathan Shewchuk for information or for pointing out information in the appendix. I also thank several readers whose comments improved this review.

\raggedright

\providecommand{\bysame}{\leavevmode\hbox to3em{\hrulefill}\thinspace}
\providecommand{\MR}{\relax\ifhmode\unskip\space\fi MR }
\providecommand{\MRhref}[2]{%
  \href{http://www.ams.org/mathscinet-getitem?mr=#1}{#2}
}
\providecommand{\href}[2]{#2}

\end{document}